\documentclass[11pt]{article}
\usepackage[english]{babel}
\usepackage{amsmath}
\usepackage{amsthm,amssymb,latexsym}
\usepackage{epic, eepic, graphics}

\newcommand\p{\circle*{0.3}}

\def\mn{{\mbox{-}}}
\def\SS{\mathcal{S}}

\newtheorem{prop}{Proposition}
\newtheorem{Lem}[prop]{Lemma}
\newtheorem{Cor}[prop]{Corollary}
\newtheorem{The}[prop]{Theorem}
\theoremstyle{definition}

\newtheorem{Exa}{Example}
\newtheorem{Rem}[prop]{Remark}
\newtheorem{Prob}{Problem}

\begin{document}
\title{Introduction to Partially Ordered Patterns}
\author{Sergey Kitaev\footnote{A part of this paper was written during the author's
stay at the Institut Mittag-Leffler,
Sweden and the University of California at San Diego.}\\
{\small Reykjav\'{\i}k University}\\ {\small Ofanleiti 2}\\
{\small IS-103 Reykjav\'{\i}k, Iceland} \\ \ \\ {\small \em
Dedicated to the 50th anniversary of my advisor, Einar
Steingr\'{\i}msson}}

\date{}
\maketitle
\thispagestyle{empty}

\begin{abstract}
 We review selected known results on partially ordered patterns (POPs) that include co-unimodal, multi- and shuffle
 patterns, peaks and valleys ((modified) maxima and minima) in permutations, the Horse permutations and others.
 We provide several new results on
 a class of POPs built on an arbitrary flat poset, obtaining, as corollaries, the bivariate generating function for the distribution of peaks (valleys)
 in permutations, links to Catalan, Narayana, and Pell numbers, as well as generalizations of a
 few results in the literature
 including the descent distribution.
 Moreover, we discuss a $q$-analogue for a result on non-overlapping segmented POPs. Finally, we suggest
 several open problems for further research.\\

{\bf Keywords:} (partially ordered) pattern, non-overlapping
occurrences, peak, valley, $q$-analogue, flat poset, co-unimodal
pattern, bijection, (exponential, bivariate) generating function,
distribution, Catalan numbers, Narayana numbers, Pell numbers
\end{abstract}

\section{Introduction and background}\label{intro}

    An occurrence of a \emph{pattern} $\tau$ in a permutation $\pi$ is defined as a subsequence
    in $\pi$ (of the same length as $\tau$) whose letters are in the same relative order as those in $\tau$.
    For example, the permutation $31425$ has three occurrences of the pattern $1\mn 2\mn 3$, namely the
    subsequences 345, 145, and 125. \emph{Generalized permutation patterns} (\emph{GPs}) being introduced
    in~\cite{BabStein} allow the requirement that some adjacent letters in a
    pattern must also be adjacent in the permutation. We indicate this requirement by removing a dash in the
    corresponding place. Say, if pattern $2\mn 31$ occurs in a
permutation $\pi$, then the letters in $\pi$
    that correspond to $3$ and $1$ are adjacent. For example, the permutation $516423$ has only one occurrence
    of the pattern $2\mn 31$, namely the subword 564, whereas the pattern $2\mn 3\mn 1$ occurs, in addition, in the
    subwords 562 and 563. Placing ``$[$'' on the left (resp.,
``$]$'' on the right) next to a pattern $p$ means the requirement
that $p$ must begin (resp., end) from the leftmost (resp.,
rightmost) letter. For example, the permutation $32415$ contains two
occurrences of the pattern $[2\mn13$, namely the subwords 324 and
315 and no occurrences of the pattern $3\mn2\mn 1]$.

\begin{figure}[h]
\begin{center}
\begin{picture}(6,2)
\put(0,0){\put(2,0){\p} \put(4,0){\p} \put(4,2){\p} \path(4,0)(4,2)

\put(1.2,0){$1'$} \put(4.5,0){1} \put(4.5,2){2}
 }
\end{picture}
\caption{A poset on three elements with the only relation $1<2$.}
\label{poset01}
\end{center}
\end{figure}
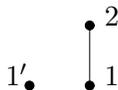

    A further generalization of the GPs is {\em partially ordered patterns} ({\em POPs}), where the letters of
    a pattern form a partially ordered set (poset), and an occurrence of such a pattern in a permutation is
    a linear extension of the corresponding poset in the order suggested by the pattern (we also pay attention
    to eventual dashes and brackets). For instance, if we have a poset on three elements labeled by $1^{\prime}$, $1$, and $2$,
    in which the only relation is $1<2$ (see Figure~\ref{poset01}), then in an occurrence of $p=1^{\prime}\mn 12$ in
    a permutation $\pi$ the letter corresponding to the $1^{\prime}$ in $p$ can be either larger or smaller than
    the letters corresponding to $12$. Thus, the permutation 31254 has three occurrences of $p$, namely $3\mn12$,
    $3\mn25$, and $1\mn25$.

    Let $\SS_n(p_1,\ldots,p_k)$ denote the set of
    $n$-permutations avoiding simultaneously each of the patterns $p_1,\ldots,p_k$.

    The POPs were introduced in~\cite{Kit2}\footnote{The POPs in this paper are the same as
    the POGPs in~\cite{Kit2}, which is
    an abbreviation for Partially Ordered Generalized Patterns.}
    as an auxiliary tool to study the maximum number of non-overlapping
    occurrences of {\em segmented} GPs ({\em SGPs}), also known as {\em consecutive} GPs, that is, the GPs, occurrences of which in permutations form contiguous subwords
    (there are no dashes).
    However, the most useful property of POPs known so far is their ability to ``encode"
    certain sets of GPs which
    provides a convenient notation for those sets and often gives an idea how to treat them.
    For example, the original
    proof of the fact that $|\SS_n(123, 132, 213)|={n\choose \lfloor n/2 \rfloor}$
    took 3 pages (\cite{Kit1}); on the other hand, if one notices that
    $|\SS_n(123, 132, 213)|=|\SS_n(11^{\prime}2)|$, where the letters $1$, $1^{\prime}$, and $2$ came from
    the same poset as above, then
    the result is easy to see. Indeed, we may use the property that the
letters in odd and even positions of a ``good" permutation do not
affect each other because of the form of $11^{\prime}2$. Thus we
choose the letters in odd positions in ${n \choose \lfloor n/2
\rfloor}$ ways, and we must arrange them in decreasing order. We
then must arrange the letters in even positions in decreasing
order too.

The POPs can be used to encode certain combinatorial objects by
restricted permutations. Examples of that are
Propositions~\ref{combobject01} and~\ref{combobject02}, as well as
several other propositions in~\cite{BurKit}. Such encoding is
interesting from the point of view of finding bijections, but it
also may have applications for enumerating certain statistics. The
idea is to encode a set of objects under consideration as a set of
permutations satisfying certain restrictions (given by certain
POPs); under appropriate encodings, this allows us to transfer the
interesting statistics from the original set to the set of
permutations, where they are easy to handle. For an illustration of
how encodings by POPs can be used, see~\cite[Thm. 2.4]{KitPM} which
deals with POPs in {\em compositions} rather than in permutations,
though, but the approach remains the same.

As a matter of fact, some POPs appeared in the literature before
they were actually introduced. Thus the notion of a POP allows us to
collect under one roof (to provide a uniform notation for) several
combinatorial structures such as {\em peaks}, {\em valleys}, {\em
modified maxima} and {\em modified minima} in permutations, {\em
Horse permutations} and $p$-{\em descents} in permutations discussed
in Section~\ref{sec1}.

This paper is organized as follows. Section~\ref{sec1} reviews
selected results in the literature related to POPs;
Section~\ref{sec2} provides a complete solution for SPOPs built on a
{\em flat poset}\footnote{The concept of a ``flat poset'' is used in
theoretical computer science~\cite{Aceto} to denote posets with one
element being less than any other element (there are no other
relations between the elements). See Figure~\ref{poset08} for the
shape of such poset.} without repeated letters. In particular, as a
corollary to a more general result, we provide the generating
function for the distribution of peaks (valleys) in permutations,
which seems to be a new result, or at least one the author could not
find in the literature (it looks like only a continued fraction
expansion of the generating function for the distribution of peaks
is known). Section~\ref{q-an} gives a $q$-analogue for a result on
non-overlapping patterns~(\cite[Thm. 16]{Kit3}). Finally, in
Section~\ref{sec3}, we state several open problems on POPs.

In what follows we need the following notations. Let $\sigma$ and
$\tau$ be two POPs of length greater than 0. We write $\sigma <\tau$
to indicate that any letter of $\sigma$ is less than any letter of
$\tau$. We write $\sigma<>\tau$ when no letter in $\sigma$ is
comparable with any letter in $\tau$. Also, {\em SPOP} abbreviates
Segmented POP.

A {\em left-to-right minimum} of a permutation $\pi$ is an element
$a_i$ such that $a_i<a_j$ for every $j<i$. Analogously we define
{\em right-to-left minimum}, {\em right-to-left maximum}, and {\em
left-to-right maximum}. If $\pi=a_1a_2\cdots a_n\in \SS_n$, then the
{\em reverse} of $\pi$ is $\pi^r:=a_n\cdots a_2a_1$, and the {\em
complement} of $\pi$ is a permutation $\pi^c$ such that
$\pi^c_{i}=n+1-a_i$, where $i\in[n]=\{1,\ldots,n\}$. We call
$\pi^r$, $\pi^c$, and $(\pi^r)^c=(\pi^c)^r$ {\em trivial
bijections}. The {\em GF} ({\em EGF}; {\em BGF}) denotes the ({\em
exponential}; {\em bivariate}) {\em generating function}.

\section{Review of selected results on POPs}\label{sec1}

In this section we review several results in the literature related
to POPs.

\subsection{Co-unimodal patterns}\label{unimod}
For a permutation $\pi=\pi_1\pi_2\cdots\pi_n\in\SS_n$, the {\em
inversion index}, $\mbox{inv}(\pi)$, is the number of ordered pairs
$(i,j)$ such that $1\leq i<j\leq n$ and $\pi_i>\pi_j$. The major
index, $\mbox{maj}(\pi)$, is the sum of all $i$ such that
$\pi_i>\pi_{i+1}$. Suppose $\sigma$ is a SPOP and
$$\mbox{place}_{\sigma}(\pi)=\{i\ |\ \pi \mbox{ has an occurrence of } \sigma \mbox{ starting at } \pi_i\}.$$
Let $\mbox{maj}_{\sigma}(\pi)$ be the sum of the elements of $\mbox{place}_{\sigma}(\pi)$.

    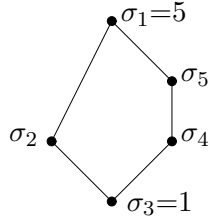
\begin{figure}[h]
\begin{center}
\begin{picture}(8,6)
\put(0,0){\put(2,2){\p} \put(4,0){\p} \put(4,6){\p} \put(6,2){\p}
\put(6,4){\p} \path(2,2)(4,0)(6,2)(6,4)(4,6)(2,2)

\put(4.3,6){$\sigma_1$=5} \put(0.6,2){$\sigma_2$}
\put(6.3,2){$\sigma_4$} \put(6.3,4){$\sigma_5$}
\put(4.5,-0.3){$\sigma_3$=1}
 }
\end{picture}
\caption{A poset for co-unimodal pattern in the case $j=3$ and
$k=5$.} \label{poset02}
\end{center}
\end{figure}
If $\sigma$ is {\em co-unimodal}, meaning that
$k=\sigma_1>\sigma_2>\cdots >\sigma_j<\cdots <\sigma_k$ for some
$2\leq j\leq k$ (see Figure~\ref{poset02} for a corresponding poset
in the case $j=3$ and $k=5$), then the following formula
holds~\cite{BjornerWachs}:
$$\sum_{\pi\in\SS_n}t^{\mbox{maj}_{\sigma}(\pi^{-1})}q^{\mbox{maj}(\pi)}=
\sum_{\pi\in\SS_n}t^{\mbox{maj}_{\sigma}(\pi^{-1})}q^{\mbox{inv}(\pi)}.$$
If $k=2$ we deal with usual descents, thus a co-unimodal pattern can
be viewed as a generalization of the notion of a descent. This may
be a reason why a co-unimodal pattern $p$ is called $p$-{\em
descent} in~\cite{BjornerWachs}. Also, setting $t=1$ we get a
well-known result by MacMahon on equidistribution of maj~and~inv.

\subsection{Peaks and valleys in permutations}

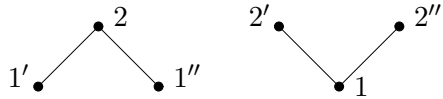
\begin{figure}[h]
\begin{center}
\begin{picture}(8,2)
\put(-4,0){\put(2,0){\p} \put(4,2){\p} \put(6,0){\p}
\path(2,0)(4,2)(6,0) \put(10,2){\p} \put(12,0){\p} \put(14,2){\p}
\path(10,2)(12,0)(14,2)

\put(1,0){$1'$} \put(6.5,0){$1''$} \put(4.5,2){2} \put(9,2){$2'$}
\put(12.5,-0.3){$1$} \put(14.5,2){$2''$}
 }
\end{picture}
\caption{Posets corresponding to peaks and valleys.} \label{poset03}
\end{center}
\end{figure}

A permutation $\pi$ has exactly $k$ {\em peaks} (resp., {\em
valleys}), also known as {\em maxima} (resp., {\em minima}), if
$|\{j\ |\ \pi_j>\max\{\pi_{j-1},\pi_{j+1}\}\}|=k$ (resp., $|\{j\ |\
\pi_j<\min\{\pi_{j-1},\pi_{j+1}\}\}|=k$). Thus, an occurrence of a
peak in a permutation is an occurrence of the SPOP $1'21''$, where
relations in the poset are $1'<2$ and $1''<2$. Similarly,
occurrences of valleys correspond to occurrences of the SPOP
$2'12''$, where $2'>1$ and $2''>1$. See Figure~\ref{poset03} for the
posets corresponding to the peaks and valleys. So, any research done
on the peak (or valley) statistics can be regarded as research on
(S)POPs (e.g., see~\cite{Warren}).

Also, results related to {\em modified maxima} and {\em modified
minima} can be viewed as results on SPOPs. For a permutation
$\sigma_1\ldots\sigma_n$ we say that $\sigma_i$ is a {\em modified
maximum} if $\sigma_{i-1}<\sigma_i>\sigma_{i+1}$ and a {\em modified
minimum} if $\sigma_{i-1}>\sigma_i<\sigma_{i+1}$, for
$i=1,\ldots,n$, where $\sigma_0=\sigma_{n+1}=0$. Indeed, we can view
a pattern $p$ as a function from the set of all symmetric groups
$\cup_{n\geq 0}\SS_n$ to the set of natural numbers such that
$p(\pi)$ is the number of occurrences of $p$ in $\pi$, where $\pi$
is a permutation. Thus, studying the distribution of modified maxima
(resp., minima) is the same as studying the function
$ab]+1'21''+[dc$ (resp., $ba]+2'12''+[cd$) where $a<b$, $c<d$ and
the other relations between the patterns' letters are taken from
Figure~\ref{poset03}. Also, recall that placing ``$[$'' (resp.,
``$]$'') next to a pattern $p$ means the requirement that $p$ must
begin (resp., end) with the leftmost (resp., rightmost) letter.

A specific result in this direction is problem 3.3.46(c) on page 195
in~\cite{GolJack}: We say that $\sigma_i$ is a {\em double rise}
(resp., {\em double fall}) if $\sigma_{i-1}<\sigma_{i}<\sigma_{i+1}$
(resp., $\sigma_{i-1}>\sigma_{i}>\sigma_{i+1}$); The number of
permutations in $\SS_n$ with $i_1$ modified minima, $i_2$ modified
maxima, $i_3$ double rises, and $i_4$ double falls is
$$\left[u_1^{i_1}u_2^{i_2-1}u_3^{i_3}u_4^{i_4}\frac{x^n}{n!}\right]\frac{e^{\alpha_2x}-e^{\alpha_1}x}
{\alpha_2e^{\alpha_1x}-\alpha_1e^{\alpha_2x}}$$ where
$\alpha_1\alpha_2=u_1u_2$, $\alpha_1+\alpha_2=u_3+u_4$.

In Corollary~\ref{valleys} we obtain explicit generating function
for the distribution of peaks (valleys) in permutations. This result
is an analogue to a result in~\cite{Entr} where the circular case of
permutations is considered, that is, when the first letter of a
permutation is though to be to the right of the last letter in the
permutation. In~\cite{Entr} it is shown that if $M(n,k)$ denotes the
number of circular permutations in $\SS_n$ having $k$ maxima, then
$$\sum_{n\geq 1}\sum_{k\geq 0}M(n,k)y^k\frac{x^n}{n!}=\frac{zx(1-z\tanh xz)}{z-\tanh xz}$$
where $z=\sqrt{1-y}$.

\subsection{Patterns containing $\Box$-symbol}

In~\cite{ManQ} the authors study simultaneous avoidance of the
patterns $1\mn3\mn 2$ and $1\Box 23$. A permutation $\pi$ avoids
$1\Box 23$ if there is no $\pi_i<\pi_j<\pi_{j+1}$ with $i<j-1$. Thus
the $\Box$ symbol has the same meaning as ``$\mn$" except for $\Box$
does not allow the letters separated by it to be adjacent in an
occurrence of the corresponding pattern. In the POP-terminology,
$1\Box 23$ is the pattern $1\mn 1'\mn 23$, or $1\mn 1'23$, or
$11'\mn 23$, where $1'$ is incomparable with the letters $1,2,$ and
$3$ which, in turn, are ordered naturally: $1<2<3$. The permutations
avoiding $1\mn3\mn 2$ and $1\Box 23$ are called {\em Horse
permutations}. The reason for the name came from the fact that these
permutations are in one to one correspondence with {\em Horse
paths}, which are the lattice paths from (0,0) to $(n,n)$ containing
the steps $(0,1)$, $(1,1)$, $(2,1)$, and $(1,2)$ and not passing the
line $y=x$. According to~\cite{ManQ}, the generating function for
the horse permutations is
$$\frac{1-x-\sqrt{1-2x-3x^2-4x^3}}{2x^2(1+x)}.$$ Moreover, in~\cite{ManQ} the generating functions for Horse permutations
avoiding, or containing (exactly) once, certain patterns are given.

In~\cite{FMan}, patterns of the form $x\mn y\Box z$ are studied,
where $xyz\in \SS_3$. Such a pattern can be written in the
POP-notation as, for example, $x\mn y\mn a\mn z$ where $a$ is not
comparable to $x$, $y$, and $z$. A bijection between permutations
avoiding the pattern $1\mn 2\Box 3$, or $2\mn 1\Box 3$, and the set
of {\em odd-dissection convex polygons} is given. Moreover,
generating functions for permutations avoiding $1\mn 3\Box 2$ and
certain additional patterns are obtained in~\cite{FMan}.

\subsection{A pattern of the form $\sigma\mn m\mn\tau$}
Let $\sigma$ and $\tau$ be two SGPs (the results below work for SPOPs as well).
We
 consider the POP $\alpha=\sigma\mn m\mn\tau$ with $m>\sigma$,
$m>\tau$, and $\sigma<>\tau$, that is, each letter of $\sigma$ is
incomparable with any letter of~$\tau$ and $m$ is the largest
letter in $\alpha$. The POP $\alpha$ is an instance of so called
{\em shuffle patterns} (see~\cite[Sec 4]{Kit2}).

\begin{The}{\rm(\cite[Thm. 16]{Kit2})}\label{shufflePatern2} Let $A(x)$, $B(x)$ and $C(x)$ be the EGF for
the number of permutations that avoid $\sigma$, $\tau$ and $\alpha$ respectively.
Then $C(x)$ is the solution to the following differential equation with $C(0)=1$:
$$C^{\prime}(x)=(A(x)+B(x))C(x) - A(x)B(x).$$
\end{The}

If $\tau$ is the empty word then $B(x)=0$ and we get the following
result for segmented GPs:

\begin{Cor}{\rm(\cite[Thm. 13]{Kit2},\cite{Knuth})}  Let $\alpha=\sigma\mn m$, where
$\sigma$ is a SGP on $[k-1]$. Let $A(x)$ {\rm(}resp., $C(x)${\rm)}
be the EGF for the number of permutations that avoid $\sigma$
{\rm(}resp., $\alpha${\rm)}. Then $C(x) = e^{F(x,A)},$ where $F(x,A)
= \int_0^x A(y)\ dy$. \end{Cor}

\begin{Exa}{\rm(\cite[Ex 15]{Kit2})} Suppose $\alpha = 12\mn 3$.
Here $\sigma = 12$, whence $A(x) = e^x$,
since there is only one permutation that avoids $\sigma$. So
$$C(x) = e^{F(x,\exp)} = e^{e^x-1}.$$ We get~\cite[Prop. 4]{Claes} since $C(x)$ is the EGF for the Bell numbers.
\end{Exa}

\begin{Cor}{\rm(\cite[Cor. 19]{Kit2})} Let $\alpha=\sigma\mn m\mn \tau$ is as described above. We consider the pattern
$\varphi (\alpha) = {\varphi}_1(\sigma)\mn m\mn {\varphi}_2(\tau)$,
where ${\varphi}_1$ and ${\varphi}_2$ are any trivial bijections.
Then $|\SS_n(\alpha)|=|\SS_n(\varphi(\alpha))|$. \end{Cor}

\subsection{Multi-patterns}

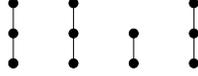
\begin{figure}[h]
\begin{center}
\begin{picture}(5,3)
\put(0,0) {\put(0,0){\p} \put(0,1){\p} \put(0,2){\p} \put(0,0){\p}
\put(2,0){\p} \put(2,1){\p} \put(2,2){\p} \put(4,0){\p}
\put(4,1){\p} \put(6,0){\p} \put(6,1){\p} \put(6,2){\p}

 \path(0,0)(0,1)(0,2) \path(2,0)(2,1)(2,2) \path(4,0)(4,1) \path(6,0)(6,1)(6,2)}
\end{picture}
\end{center}
\caption{A poset corresponding to a multi-pattern.} \label{poset09}
\end{figure}

Suppose $\{ {\sigma}_1, {\sigma}_2, \ldots , {\sigma}_k \}$ is a set
of segmented GPs and $p={\sigma}_1\mn {\sigma}_2\mn \cdots\mn
{\sigma}_k$ where each letter of ${\sigma}_i$ is incomparable with
any letter of ${\sigma}_j$ whenever $i \neq j$
(${\sigma}_i<>{\sigma}_j$). We call such POPs {\em multi-patterns}.
Clearly, the Hasse diagram for such a pattern is $k$ disjoint chains
similar to that in Figure~\ref{poset09}.

\begin{The}{\rm(\cite[Thm. 23 and Cor. 24]{Kit2})}\label{mult1} The number of permutations avoiding the pattern
$p={\sigma}_1\mn {\sigma}_2\mn \cdots\mn {\sigma}_k$ is equal to
that avoiding a multi-pattern obtained from $p$ by an arbitrary
permutation of $\sigma_i$'s as well as by applying to $\sigma_i$'s
any of trivial bijections. \end{The}

The following theorem is the basis for calculating the number of
permutations that avoid a multi-pattern.

\begin{The}{\rm(\cite[Thm. 28]{Kit2})}\label{mult2} Let
$p={\sigma}_1\mn {\sigma}_2\mn \cdots\mn {\sigma}_k$ be a
multi-pattern and let $A_i(x)$ be the EGF for the number of
permutations that avoid ${\sigma}_i$. Then the EGF $A(x)$ for the
number of permutations that avoid~$p$ is
$$A(x)=\displaystyle\sum_{i=1}^{k}A_i(x)\displaystyle\prod_{j=1}^{i-1}((x-1)A_j(x) + 1).$$ \end{The}

\begin{Cor}{\rm(\cite[Cor. 26]{Kit2})} Let $p={\sigma}_1\mn {\sigma}_2\mn \cdots\mn {\sigma}_k$ be a multi-pattern, where $|{\sigma}_i|=2$ for all $i$.
That is, each ${\sigma}_i$ is either 12 or 21. Then the EGF for the
number of permutations that avoid $p$ is given by
$$A(x)=\frac{1-(1+(x-1)e^x)^k}{1-x}.$$ \end{Cor}

\begin{Rem}\label{goodrem} Although the results in Theorems~\ref{mult1} and~\ref{mult2} are stated in~\cite{Kit2} for
$\sigma_i$'s which are SGPs, one can see that the same arguments
work for $\sigma_i$'s which are SPOPs. Thus we have a generalization
of these theorems. \end{Rem}

\subsection{Non-overlapping patterns -- an application of POPs}
This subsection deals additionally with occurrences of patterns in
words. The letters $1,2,1',2'$ appearing in the examples below are
ordered as in Figure~\ref{poset001} to be found on
page~\pageref{poset001}.

Theorem~\ref{mult2} and its counterpart in the case of
words~\cite[Thm. 4.3]{KitMan} and~\cite[Cor. 4.4]{KitMan}, as well
as Remark~\ref{goodrem} applied for these results, give an
interesting application of the multi-patterns in finding a certain
statistic, namely the {\em maximum number of non-overlapping
occurrences of a SPOP} in permutations and words. For instance, the
maximum number of non-overlapping occurrences of the SPOP $11'2$ in
the permutation 621394785 is 2, and this is given by the occurrences
$213$ and $478$, or the occurrences $139$ and $478$.

Theorem~\ref{dist} generalizes~\cite[Thm. 32]{Kit2} and~\cite[Thm.
5.1]{KitMan}.

\begin{The}{\rm(\cite[Thm. 16]{Kit3})}\label{dist} Let $p$ be a SPOP and $B(x)$ {\rm(}resp.,
$B(x;k)${\rm)} is the EGF {\rm(}resp., GF{\rm)} for the number of
permutations {\rm(}resp., words over $[k]${\rm)} avoiding $p$. Let
$D(x,y)=\sum_{\pi}y^{N(\pi)}\frac{x^{|\pi|}}{|\pi|!}$ and
$D(x,y;k)=\sum_{n\geq 0}\sum_{w\in [k]^n}y^{N(w)}x^n$ where $N(s)$
is the maximum number of non-overlapping occurrences of $p$ in~$s$.
Then $D(x,y)$ and $D(x,y;k)$ are given by
$$\frac{B(x)}{1-y(1+(x-1)B(x))} \mbox{\ \ \ and
\ \ \ } \frac{B(x;k)}{1-y(1+(kx-1)B(x;k))}.$$
\end{The}
The following examples are corollaries to Theorem~\ref{dist}.

\begin{Exa}{\rm(\cite[Ex 1]{Kit3})} If we consider the SPOP $11'$ then clearly $B(x)=1+x$
and $B(x;k)=1+kx$. Hence, $$D(x,y)=\frac{1+x}{1-yx^2}=\sum_{i\geq
0}(x^{2i} + x^{2i+1})y^i,$$and
$$D(x,y;k)=\frac{1+kx}{1-y(kx)^2}=\sum_{i\geq 0}((kx)^{2i} +
(kx)^{2i+1})y^i.$$
\end{Exa}

\begin{Exa}{\rm(\cite[Ex 2]{Kit3})} For permutations, the distribution
of the maximum number of non-overlapping occurrences of the SPOP
$122'1'$ is given by
$$D(x,y)=\frac{\frac{1}{2}+\frac{1}{4}\tan x(1+e^{2x}+2e^{x}\sin
x)+\frac{1}{2}e^x\cos x}{1-y(1+(x-1)(\frac{1}{2}+\frac{1}{4}\tan
x(1+e^{2x}+2e^{x}\sin x)+\frac{1}{2}e^x\cos x))}.$$
\end{Exa}

\subsection{Segmented patterns of length four}
In this subsection we provide the known results related to SPOPs of
length four. Corollaries~\ref{newfour1} and~\ref{newfour} in
subsection~\ref{subsec1} give extra results in this direction. In
subsection~\ref{usolved-4} we provide unsolved cases with initial
values for the number of the restricted permutations. In this
subsection, $A(x)=\sum_{n\geq 0}A_nx^n/n!$
 is the EGF for the number of permutations in question. The patterns
in the subsection are built on the poset from Figure~\ref{poset001}
and the letter $1''$ is not comparable to any other letter.

\begin{The}{\rm(\cite[Thm. 30]{Kit2})} For the SPOP $122^{\prime}1^{\prime}$, we have that
$$A(x) = \frac{1}{2} + \frac{1}{4}\tan x (1+e^{2x} + 2e^{x}\sin x) + \frac{1}{2}e^{x} \cos x .$$
\end{The}

\begin{prop}{\rm(\cite[Prop. 8,9]{Kit3})}\label{combobject01} There are ${n-1 \choose \lfloor (n-1)/2 \rfloor}{n \choose \lfloor n/2 \rfloor}$
permutations in $\SS_n$ that avoid the SPOP
$12^{\prime}21^{\prime}$. The $(n+1)$-permutations avoiding $12'21'$
are in one-to-one correspondence with different walks of $n$ steps
between lattice points, each in a direction N, S, E or W, starting
from the origin and remaining in the positive quadrant.
\end{prop}

\begin{prop}{\rm(\cite[Prop. 4,5,6]{Kit3})} For the SPOP $11'1''2$, one has
$$A_n=\frac{n!}{\lfloor n/3 \rfloor !\lfloor (n+1)/3 \rfloor !\lfloor (n+2)/3 \rfloor !},$$ and
for the SPOP $11'21''$ and $n\geq 1$, we have
$A_n=n\cdot\displaystyle{n-1 \choose \lfloor (n-1)/2 \rfloor}$.
Moreover, for the SPOPs $1'1''12$ and $1'121''$, we have
$A_0=A_1=1$, and, for $n\geq 2$, $A_n=n(n-1)$.\end{prop}

\begin{prop}{\rm(\cite[Prop. 7]{Kit3})} For the SPOP $1231'$, we have
$$A(x)=xe^{x/2}\left(\cos \frac{\sqrt{3}x}{2}
-\frac{\sqrt{3}}{3}\sin \frac{\sqrt{3}x}{2} \right)^{-1} + 1,$$ and
for the SPOPs $1321'$ and $2131'$, we have
$$A(x)=x(1-\int_{0}^{x}e^{-t^2/2}\ dt)^{-1} +1.$$\end{prop}

We end up this subsection with a result on multi-avoidance of SPOPs
that has a combinatorial interpretation.

\begin{prop}{\rm(\cite[Prop. 2.1,2.2]{BurKit})}\label{combobject02}
There are $2\binom{n}{\lfloor n/2 \rfloor}$ $n$-permutations
avoiding the SPOPs $11'22'$ and $22'11'$ simultaneously. For $n\ge
3$, there is a bijection between such $n$-permutations and the set
of all $(n+1)$-step walks on the $x$-axis with the steps $a=(1,0)$
and $\bar{a}=(-1,0)$ starting from the origin but not returning to
it.
\end{prop}

\section{Patterns built on flat posets}\label{sec2}

\begin{figure}[h]
\begin{center}
\begin{picture}(5,3)
\put(0,0) {\put(2,0){\p} \put(0,2){\p} \put(1,2){\p} \put(2,2){\p}
 \put(3,2){\p} \put(4,2){\p} \put(2,-0.6){$a$} \put(-0.5,2.5){$a_1$}
 \put(0.5,2.5){$a_2$} \put(2,2.5){$\cdots$}  \put(4,2.5){$a_k$}
 \path(0,2)(2,0)(1,2) \path(2,2)(2,0)(3,2) \path(2,0)(4,2)}
\end{picture}
\end{center}
\caption{A flat poset.} \label{poset08}
\end{figure}
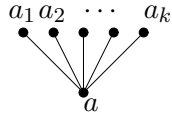

In this section, we consider flat posets built on $k+1$ elements
$a,a_1,\ldots, a_k$ with the only relations $a<a_i$ for all $i$. A
Hasse diagram for the flat poset is in Figure~\ref{poset08}.

\subsection{Avoidance and distribution of the patterns}\label{subsec1}

The following proposition generalizes~\cite[Prop. 6]{Claes}. Indeed,
letting $k=2$ in the proposition we deal with involutions and
permutations avoiding $1\mn 23$ and $1\mn 32$.

\begin{prop}\label{claes} The permutations in ${\SS}_n$ having cycles of length at most $k$ are in
one-to-one correspondence with permutations in ${\SS}_n$ that avoid
$a\mn a_1\cdots a_k$. \end{prop}

\begin{proof} We construct a bijection in a similar to~\cite[Prop. 6]{Claes} way.

Let $\pi\in {\SS}_n$ be a permutation with cycles of length at most
$k$. A standard form for writing $\pi$ in cycle notation is
requiring that
\begin{itemize}
\item[(1)] Each cycle is written with its least element first;
\item[(2)] The cycles are written in decreasing order of their least element.
\end{itemize}
Define $\hat{\pi}$ to be the permutation obtained from $\pi$ by
writing it in standard form and erasing the parentheses separating
the cycles. The permutation $\hat{\pi}$ avoids $a\mn a_1\cdots a_k$.
Indeed, the distance between two left-to-right minima (the number of
letters between them) in $\hat{\pi}$ does not exceed $k-1$ because
of the restriction on the cycle lengths. Thus if $\hat{\pi}$
contains $a\mn a_1\cdots a_k$ then among the letters of $\hat{\pi}$
corresponding to $a_1\cdots a_k$ there is at least one left-to-right
minimum, say $m$, and the letter in $\hat{\pi}$ corresponding to $a$
must be less than $m$. This contradicts the definition of a
left-to-right minimum.

Conversely, if $\hat{\pi}$ is an $a\mn a_1\cdots a_k$-avoiding
permutation then any two of its consecutive left-to-right minima are
at distance not exceeding $k-1$ from each other, since otherwise we
have an occurrence of $a\mn a_1\cdots a_k$ starting at a
left-to-right minimum preceding a factor of length at least $k$ that
does not contain other left-to-right minima. The left-to-right
minima of $\hat{\pi}$ define cycles of~$\pi$.
\end{proof}

\begin{Cor}\label{nicecor} The EGF for the number of permutations avoiding $a\mn a_1\cdots a_k$ is given by
$\exp(\sum_{i=1}^{k}x^i/i)$.\end{Cor}

\begin{proof} According to Proposition~\ref{claes} we only need to find the EGF $p(x)=\sum_{n\geq 0}p_nx^{n}/n!$
for the number of permutations with cycles of length at most $k$,
which is known (see, e.g.,~\cite{GesStan}), but we rederive it here.

Suppose $\pi$ is an $n$-permutation with cycles of length at most
$k$ and 1 occurs in a cycle $C$. If $i$ is the number of neighbors
of 1 in $C$ then $0\leq i\leq k-1$ and there are ${n-1\choose i}i!$
possibilities for choosing such $C$. Thus
$$p_n=\displaystyle\sum_{i=0}^{k-1}{n-1\choose i}i!p_{n-i-1}$$ which
after summing over all $n\geq 1$ gives
$$p'(x)=(1+x+\cdots+x^{k-1})p(x)$$ and therefore the claim is true
since $P(0)=1$.
\end{proof}

\begin{prop}\label{easyprop1}
One has ${\SS}_n(a\mn a_1\cdots a_k)={\SS}_n(aa_1\cdots a_k)$, and
thus the EGF for the number of permutations avoiding $aa_1\cdots
a_k$ is $\exp(\sum_{i=1}^{k}x^i/i)$.
\end{prop}

\begin{proof} Clearly ${\SS}_n(a\mn a_1\cdots a_k)\subseteq {\SS}_n(aa_1\cdots a_k)$. Suppose now that $\pi\in {\SS}_n(aa_1\cdots a_k)$ and $\pi$ contains an occurrence of
$a\mn a_1\cdots a_k$, say $\pi_i\pi_j\pi_{j+1}\cdots\pi_{j+k-1}$
where $i+1<j$. We will get a contradiction which will show that
${\SS}_n(aa_1\cdots a_k)\subseteq {\SS}_n(a\mn a_1\cdots a_k)$.

One can assume that $j-i$ is minimal out of all occurrences of $a\mn
a_1\cdots a_k$ in $\pi$. If $\pi_{j-1}<\pi_i$ then
$\pi_{j-1}\pi_j\pi_{j+1}\cdots\pi_{j+k-1}$ is an occurrence of
$aa_1\cdots a_k$, a contradiction to $\pi\in {\SS}_n(aa_1\cdots
a_k)$; otherwise, $\pi_i\pi_{j-1}\pi_j\cdots\pi_{j+k-2}$ is an
occurrence of $a\mn a_1\cdots a_k$, a contradiction to $j-i$ being
minimal.
\end{proof}

\begin{Cor}\label{newfour1} The EGF for the number of permutations avoiding $aa_1a_2a_3$ is given by
$\exp(x+x^2/2+x^3/3)$.\end{Cor}

\begin{The}{\rm(Distribution of $aa_1a_2\cdots a_k$)}\label{spogpdist} Let
$$P:=P(x,y)=\sum_{n\geq 0}\sum_{\pi\in\SS_n}y^{e(\pi)}x^{n}/n!$$ be
the BGF on permutations, where $e(\pi)$ is the number of occurrences
of the SPOP $p=aa_1a_2\cdots a_k$ in $\pi$. Then $P$ is the solution
of
\begin{equation}\label{eq1} \frac{\partial P}{\partial x}=yP^2+\frac{(1-y)(1-x^k)}{1-x}P\end{equation}
with the initial condition $P(0,y)=1$. \end{The}

\begin{proof} Suppose $\pi=\pi'1\pi''$ is a permutation. Then
\[ e(\pi) = \left\{ \begin{array}{ll} e(\pi')+e(\pi'')+1 & \mbox{if $|\pi''|\geq k$,}
\\ e(\pi') & \mbox{if $|\pi''|< k$}
\end{array}
\right. \] since an occurrence of $p$ cannot start at $\pi'$ and end
not in $\pi'$; also when $\pi''$ is of length at least $k$ it
contributes one extra occurrence of $p$ starting at 1.

Suppose
$P_{<k}:=P_{<k}(x,y)=\sum_{n=0}^{k-1}\sum_{\pi\in\SS_n}y^{e(\pi)}x^{n}/n!=\sum_{n\geq
0}x^n=\frac{1-x^k}{1-x}$. Readers familiar with the symbolic method
can now see that
$$P'=P(y(P-P_{<k})+P_{<k})$$
with the initial condition $P(0,y)=1$ and the desired is easy to get
by plugging in $P_{<k}$ and rewriting the equation.

The rest of the proof is dedicated to a brief explanation of the
symbolic method (see~\cite{FlaSed} for more details) and applying it
to our case. In our presentation we follow~\cite{ElizNoy}.

There is a direct correspondence between set-theoretic operations on
combinatorial classes and algebraic operations on EGFs. Let
$\mathcal{A}$, $\mathcal{B}$, and $\mathcal{C}$ be classes of
labeled combinatorial objects, and $A(x)$, $B(x)$, and $C(x)$ be
their EGFs respectively. Then if $\mathcal{A}=\mathcal{B}\cup
\mathcal{C}$ is the union of disjoint copies then $A(x)=B(x)+C(x)$;
if $\mathcal{A}=\mathcal{B}\star \mathcal{C}$ is the labeled
product, that is, the usual Cartesian product enriched with the
relabeling operation, then $A(x)=B(x)C(x)$; if
$\mathcal{A}=\mathcal{B}^{\Box}\star \mathcal{C}$ is the box
product, that is, the subset of $\mathcal{B}\star \mathcal{C}$
formed by those pairs in which the smallest label lies in the
$\mathcal{B}$ component, then $A(x)=\int_0^x(\frac{d}{dt}B(t))\cdot
C(t)dt$. The same holds if we have the BGFs instead of EGFs.

Let $\mathcal{P}$ be the class of all permutations and
$\mathcal{P}_{<k}$ is the class of permutations of length less than
$k$. With some abuse of notation, we introduce the parameter $y$ in
the equation for classes meaning that it will be placed there when
we write the corresponding differential equations for the BGFs. With
this notation and using the property of $e(\pi)$, we can write
$$\mathcal{P}=\{\epsilon\}+\{x\}^{\Box}\star\mathcal{P}\star[y(\mathcal{P}-\mathcal{P}_{<k})+\mathcal{P}_{<k}]$$
where $\epsilon$ is the empty permutation. We differentiate the
corresponding equation for BFGs to get the desired result.
\end{proof}

Note, that if $y=0$ in Theorem~\ref{spogpdist} then the function in
Corollary~\ref{nicecor}, due to Proposition~\ref{easyprop1}, is
supposed to be the solution to~(\ref{eq1}), which is true. If $k=1$
in Theorem~\ref{spogpdist}, then as the solution to~(\ref{eq1}) we
get nothing else but the distribution of descents in permutations:
$(1-y)(e^{(y-1)x}-y)^{-1}$. Thus Theorem~\ref{spogpdist} can be
thought as a generalization of the result on the descent
distribution.

The following theorem generalizes Theorem~\ref{spogpdist}. Indeed,
Theorem~\ref{spogpdist} is obtained from Theorem~\ref{spogpdist1} by
plugging in $\ell=0$ and observing that obviously $aa_1\cdots a_k$
and $a_1\cdots a_ka$ are equidistributed.

\begin{The}{\rm(Distribution of $a_1a_2\cdots a_kaa_{k+1}a_{k+2}\cdots a_{k+\ell}$)}\label{spogpdist1} Let
$$P:=P(x,y)=\sum_{n\geq 0}\sum_{\pi\in\SS_n}y^{e(\pi)}x^{n}/n!$$ be
the BGF of permutations where $e(\pi)$ is the number of occurrences
of the SPOP $p=a_1a_2\cdots a_kaa_{k+1}a_{k+2}\cdots a_{k+\ell}$ in
$\pi$. Then $P$ is the solution of
\begin{equation}\label{eq2}
\frac{\partial P}{\partial
x}=y\left(P-\frac{1-x^k}{1-x}\right)\left(P-\frac{1-x^{\ell}}{1-x}\right)+\frac{2-x^k-x^{\ell}}{1-x}P-\frac{1-x^k-x^{\ell}+x^{k+\ell}}{(1-x)^2}.\end{equation}
with the initial condition $P(0,y)=1$. \end{The}

\begin{proof} A proof is straightforward applying the technique introduced in the proof of
Theorem~\ref{spogpdist}. We use the same notation and adjusted steps
of that proof without explanations.

Suppose $\pi=\pi'1\pi''$ is a permutation. Then
\[ e(\pi) = \left\{ \begin{array}{ll} e(\pi')+e(\pi'')+1 & \mbox{if $|\pi'|\geq k$ and $|\pi''|\geq \ell$ ,}
\\ e(\pi') + e(\pi'') & \mbox{otherwise.}
\end{array}
\right. \]

One can now see that $\mathcal{P}$ is equal to
$$\{\epsilon\}+\{x\}^{\Box}\star[y(\mathcal{P}-\mathcal{P}_{<k})\star(\mathcal{P}-\mathcal{P}_{<\ell})
+(P-\mathcal{P}_{<k})\star \mathcal{P}_{<\ell} +
\mathcal{P}_{<k}\star (\mathcal{P}-\mathcal{P}_{<\ell})+
\mathcal{P}_{<k}\star\mathcal{P}_{<\ell}]$$ and the rest is obtained
by rewriting in terms of BGFs and differentiating.
\end{proof}

If $y=0$ in Theorem~\ref{spogpdist1} then we get the following
corollary:

\begin{Cor}\label{generate1} The EGF $A(x)=\sum_{n\geq 0}A_nx^n/n!$ for the number of permutations avoiding the
SPOP $p=a_1a_2\cdots a_kaa_{k+1}a_{k+2}\cdots a_{k+\ell}$ satisfies
the following differential equation with the initial condition
$A(0)=1$:
$$A'(x)=\frac{2-x^k-x^{\ell}}{1-x}A(x)-\frac{1-x^k-x^{\ell}+x^{k+\ell}}{(1-x)^2}.$$
\end{Cor}

The following corollaries to Corollary~\ref{generate1} are obtained
by plugging in $k=\ell=1$ and $k=1$ and $\ell=2$ respectively.

\begin{Cor}\rm{(}\cite{Kit1}\rm{)}\label{valleyless} The EGF for the number of permutations avoiding $a_1aa_2$ is
$(\exp(2x)+1)/2$ and thus $|\SS_n(a_1aa_2)|=2^{n-1}$.\end{Cor}

\begin{Cor}\label{newfour} The EGF for the number of permutations avoiding $a_1aa_2a_3$ is
$$1+\sqrt{\frac{\pi}{2}}(\mbox{erf}(\frac{1}{\sqrt{2}}x+\sqrt{2})-\mbox{erf}(\sqrt{2}))e^{\frac{1}{2}x(x+4)+2}$$
where $\mbox{erf}(x)=\frac{2}{\sqrt{\pi}}\displaystyle\int_{0}^{x}e^{-t^2}\ dt$ is the error
function.
\end{Cor}

If $k=1$ and $\ell=1$ then our pattern $a_1aa_2$ is nothing else but
the valley statistic. In~\cite{RieZel} a recursive formula for the
generating function of permutations with exactly $k$ valleys is
obtained, which however does not seem to allow (at least easily)
finding the corresponding BGF. As a corollary to
Theorem~\ref{spogpdist1} we get the following BGF by
solving~(\ref{eq2}) for $k=1$ and $\ell=1$:

\begin{Cor}\label{valleys} The BGF for the distribution of peaks (valleys)
in permutations is given by
$$1-\frac{1}{y}+\frac{1}{y}\sqrt{y-1}\cdot \tan\left(x\sqrt{y-1}+\arctan\left(\frac{1}{\sqrt{y-1}}\right)\right).$$\end{Cor}

Expanding the BGF in Corollary~\ref{valleys} we can get, for
example, the sequences A000431, A000487, and A000517 appearing
in~\cite{Sloane} for the number of permutations with exactly one,
two, and three valleys respectively. Note, that we have already
obtained the number of {\em valleyless} permutations in
Corollary~\ref{valleyless}. The valleyless permutations were studied
in~\cite{RieZel}.

\subsection{Distribution of the patterns with additional restrictions}

The results from this subsection are in a similar direction as that
in the papers~\cite{BraClaSte},~\cite{Man}, \cite{ManVan}, and
several other papers, where the authors study $1\mn 3\mn 2$-avoiding
permutations with respect to avoidance/count of other patterns. Such
a study not only gives interesting enumerative results, but also
provides a number of applications (see~\cite{BraClaSte}).

To state the theorem below, we define
$P_{k}=\sum_{n=0}^{k-1}\frac{1}{n+1}{2n \choose n}x^n$. That is,
$P_{k}$ is the sum of initial $k$ terms in the expansion of the
generating function $\frac{1-\sqrt{1-4x}}{2x}$ of the Catalan
numbers.

\begin{The} {\rm(Distribution of $a_1a_2\cdots a_kaa_{k+1}a_{k+2}\cdots a_{k+\ell}$ on $\SS_n(2\mn 1\mn 3)$)}\label{spogpdist2} Let
$$P:=P(x,y)=\sum_{n\geq 0}\ \sum_{\pi\in\SS_n(2\mn 1\mn 3)}y^{e(\pi)}x^{n}$$ be
the BGF of $2\mn 1\mn 3$-avoiding permutations where $e(\pi)$ is the
number of occurrences of the SPOP $p=a_1a_2\cdots
a_kaa_{k+1}a_{k+2}\cdots a_{k+\ell}$ in $\pi$. Then $P$ is given by
$$\frac{1-x(1-y)(P_{k}+P_{\ell})-\sqrt{(x(1-y)(P_{k}+P_{\ell})-1)^2-4xy(x(y-1)P_{k}P_{\ell}+1)}}{2xy}.$$\end{The}

\begin{proof} Let $\pi=\pi_11\pi_2\in\SS_n(2\mn 1\mn 3)$. Then each letter in $\pi_1$ must be greater than any
letter in $\pi_2$, where both $\pi_1$ and $\pi_2$ must necessarily
be $2\mn 1\mn 3$-avoiding. Conversely, every permutation of this
form is clearly $2\mn 1\mn 3$-avoiding.

It is easy to see that
$e(\pi)=e(\pi_1)+e(\pi_2)+\delta_{|\pi_1|,|\pi_2|}$, where
\[ \delta_{|\pi_1|,|\pi_2|} = \left\{ \begin{array}{ll} 1 & \mbox{if $|\pi_1|\geq k$ and $|\pi_2|\geq \ell$ ,}
\\ 0 & \mbox{otherwise.}
\end{array}
\right. \] Using the symbolic method we get that, in terms of GFs,
$$P=1+x(y(P-P_k)(P-P_{\ell})+P_k\cdot P+P\cdot P_{\ell}-P_k\cdot P_{\ell})$$
where 1 corresponds to the empty permutation, and we subtracted
$P_k\cdot P_{\ell}$ since the permutations corresponding to this
term are counted twice, namely in $P_k\cdot P$ and in $P\cdot
P_{\ell}$.

To get the desired we solve the equation above for $P$.
\end{proof}

We now discuss several corollaries to Theorem~\ref{spogpdist2}. Note
that letting $y=1$ we obtain the GF for the Catalan numbers. Also,
letting $y=0$ in the expansion of $P$, we obtain the GF for the
number of permutations avoiding simultaneously the patterns $2\mn
1\mn 3$ and $a_1a_2\cdots a_kaa_{k+1}a_{k+2}\cdots a_{k+\ell}$.

If $k=1$ and $\ell=0$ in Theorem~\ref{spogpdist2}, then $P_k=1$ and
$P_{\ell}=0$, and we obtain the distribution of descents in $2\mn
1\mn 3$-avoiding permutations. This distribution gives the {\em
triangle of Narayana numbers} (see~\cite[A001263]{Sloane} for more
details).

If $k=\ell=1$ in Theorem~\ref{spogpdist2} then we deal with avoiding
the pattern $2\mn 1\mn 3$ and counting occurrences of the pattern
$312$, since any occurrence of $a_1aa_2$ in a legal permutation must
be an occurrence of $312$ and vice versa. Thus the BGF of $2\mn 1\mn
3$-avoiding permutations with a prescribed number of occurrences of
$312$ is given by $$\frac{1-2x(1-y)-\sqrt{4(1-y)x^2+1-4x}}{2xy}.$$

Reading off the coefficients of the terms involving only $x$ in the
expansion of the function above, we can see that the number of
$n$-permutations avoiding simultaneously the patterns $2\mn 1\mn 3$
and $312$ is $2^{n-1}$, which is known and is easy to see directly
from the structure of such permutations.

Reading off the coefficients of the terms involving $y$ to the power
1 we see that the number of $n$-permutations avoiding $2\mn 1\mn 3$
and having exactly one occurrence of the pattern $312$ is given by
$(n-1)(n-2)2^{n-4}$. The corresponding sequence appears
as~\cite[A001788]{Sloane} and gives an interesting fact which we
state as Proposition~\ref{faces}. We give a combinatorial proof of
that fact.

\begin{prop}\label{faces} There is a bijection between 2-dimensional faces in
the $(n+1)$-dimensional hypercube and $2\mn 1\mn 3$-avoiding
$(n+2)$-permutations with exactly one occurrence of the pattern
$312$. \end{prop}

\begin{proof} Recall that a node in a hypercube is at level $i$ if the binary vector corresponding to it
contains $i$ 1's.

A 2-dimensional face in $(n+1)$-dimensional hypercube can be
specified by choosing two positions in an $(n+1)$-binary vector and
fixing the remaining entries of the vector to be 0 or 1 (in
$2^{n-1}$ ways). Indeed, any 2-dimensional face in a hypercube is a
4-cycle having two nodes at the same, say $i$-th, level, one node at
the $(i+1)$-st level and one node at the $(i-1)$-st level. Moreover,
the binary vectors corresponding to the nodes from the $i$-th level
must differ only in two coordinates and thus one of the vectors has
1 and 0 in these coordinates whereas the second vector has 0 and 1
there. So, the number of 2-dimensional faces in the
$(n+1)$-dimensional hypercube is given by ${n+1 \choose 2}2^{n-1}$
which is the same as the number of the $(n+2)$-permutations under
consideration (we refer to such permutations as ``good
permutations").

We now describe the structure of the good permutations. Suppose
$\pi=\pi_11\pi_2$ is a good permutation. Clearly, to avoid $2\mn
1\mn 3$, any letter of $\pi_1$ must be greater than any letter of
$\pi_2$. If $\pi_1$ and $\pi_2$ are non-empty, then the unique
occurrence of the pattern $312$ involves 1 in $\pi$ and both $\pi_1$
and $\pi_2$ must avoid simultaneously $2\mn 1\mn 3$ and $312$. The
permutations avoiding both $2\mn 1\mn 3$ and $312$ have one peak,
that is, the elements to the right (resp., left) of the largest
element must be in decreasing (resp., increasing) order. If $\pi_1$
(resp., $\pi_2$) is empty, then $\pi_2$ (resp., $\pi_1$) is a good
permutation and we use induction on length to describe the structure
of $\pi$.

Given a 2-dimensional face defined by a binary vector $${\bf
v}=a_1\cdots a_ixa_{i+2}\cdots a_jya_{j+2}\cdots a_{n+1}$$
 with chosen positions $i+1$ and $j+1$ filled by $x$ and $y$
(if $y$ is next to $x$ then $j=i+1$; if $x$ is the leftmost element
then $i=0$; if $y$ is the rightmost element then $j=n$). Based on
the structure considerations above, we describe a procedure to find
a good $(n+2)$-permutation corresponding to ${\bf v}$. We read ${\bf
v}$ from left to right and place $1,2,\ldots, n+2$, one by one, into
our permutation $\pi=\pi_1\pi_2\cdots\pi_{n+2}$ which we think of as
being initially $n+2$ empty slots. If we write, say, $\pi'_k$ then
we mean that the $k$-th slot of $\pi$ is filled.

We start filling $\pi$ by reading $a_k$, $k=1,2,\ldots,i$: if
$a_k=0$, place $k$ into the leftmost empty slot of $\pi$; place $k$
into the rightmost empty slot otherwise. Suppose that as the result
of filling the first $i$ elements we get
$\pi'_1\cdots\pi'_t\pi_{t+1}\cdots
\pi_{n+t-i+2}\pi'_{n+t-i+3}\cdots\pi'_{n+2}$. Set
$\pi_{t+j-i+1}=i+1$. Note that currently we have the word
$\pi'_1\cdots\pi'_tA(i+1)B\pi'_{n+t-i+3}\cdots\pi'_{n+2}$, where $A$
and $B$ consist of empty slots, $|A|=j-i\geq 1$ and $|B|=n-j+1\geq
1$. In what follows, any element to be filled in $A$ is greater than
any element to be filled in $B$, and thus the element $i+1$ is
involved in an occurrence of the pattern $312$. This occurrence will
be the only one in the permutation.

We fill in $B$ by reading $a_k$, $k=j+2,\ldots, n+1$ and placing the
elements $(i+2),\ldots,(n-j+i+1)$, one by one, as follows: if
$a_k=0$, place the current element into the leftmost empty slot of
$B$; place the current element insto the rightmost empty slot
otherwise. We place $(n-j+i+2)$ in the remaining empty slot of $B$.
Fill in the remaining elements, one by one in increasing order, into
$A$ by reading $a_k$, $k=i+2,\ldots,j$ in the way similar to that
when proceeding with $B$. In particular, $n+2$ will be placed in the
remaining empty slot of $A$.

For example, the face $110x0y01$ corresponds to the permutation
$389457621$, where $A$ is filled by 89 and $B$ by 576.

Our map is obviously injective and the converse to it is easy to
see.
\end{proof}

If $k=1$ and $\ell=2$ in Theorem~\ref{spogpdist2} then we deal with
avoiding the pattern $2\mn 1\mn 3$ and counting occurrences of the
pattern $a_1aa_2a_3$. In particular, one can see that the number of
permutations avoiding simultaneously $2\mn 1\mn 3$ and $a_1aa_2a_3$
is given by the {\em Pell numbers} $p(n)$ defined as $p(n) = 2p(n-1)
+ p(n-2)$ for $n>1$; $p(0) = 0$ and $p(1) = 1$. The Pell numbers
appear as~\cite[A000129]{Sloane}, where one can find objects related
to our restricted permutations.

\section{$q$-analogues for non-overlapping SPOPs}\label{q-an}

The purpose of this section is to prove Theorem~\ref{q-analog} which
is a $q$-analogue of~\cite[Thm. 16]{Kit3}. In fact, the formulation
of Theorem~\ref{q-analog} is similar to that of the $q$-analogue
of~\cite[Thm. 32]{Kit2} obtained in~\cite{MenRem}. Moreover, to
prove Theorem~\ref{q-analog} one can use the same arguments as those
in~\cite{MenRem} involving rather complicated considerations based
on symmetric functions, but we choose a simpler proof that is
similar to proving~\cite[Thm. 32]{Kit2} in~\cite{Kit2}.

We fix some notations. Let $p$ be a segmented POP (SPOP) and
$A^p_{n,k}$ be the number of $n$-permutations avoiding $p$ and
having $k$ inversions. As usually, $[n]_q=q^0+\cdots +q^{n-1}$,
$[n]_q!=[n]_q\cdots [1]_q$, $\left[\begin{array}{c} n \\ i
\end{array}\right]_q= \frac{[n]_q!}{[i]_q![n-i]_q!}$, and, as above,
$\mbox{inv}(\pi)$ denotes the number of inversions in a permutation
$\pi$. We set $A^p_n(q)=\sum_{\pi \mbox{ avoids }
p}q^{\mbox{inv}(\pi)}$. Moreover,
$$A^p_q(x)=\sum_{n,k}A^p_{n,k}q^k\frac{x^n}{[n]_q!}=\sum_{n}A^p_n(q)\frac{x^n}{[n]_q!}=
\sum_{\pi \mbox{ avoids }p}q^{\mbox{inv}(\pi)}
\frac{x^{|\pi|}}{[|\pi|]_q!}.$$

All the definitions above are similar in case of permutations that
{\em quasi-avoid} $p$, indicated by $B$ rather than $A$, namely,
those permutations that have exactly one occurrence of $p$ and this
occurrence consists of the $|p|$ rightmost letters in the
permutations.

\begin{Lem}\label{lem01}\rm{(A $q$-analogue of~\cite[Prop. 4]{Kit2} that is valid for POPs)}
We have $B^p_q(x)=(x-1)A^p_q(x)+1$. \end{Lem}

\begin{proof} If we consider all $(n-1)$-permutations avoiding $p$
(the number of those, if we register inversions, is $A^p_{n-1}(q)$)
and all possible extensions of these permutations to the
$n$-permutations by writing one more letter to the right; then the
number of obtained permutations, with inversions registered, is
$(1+q+\cdots+q^{n-1})A_{n-1}(q)=[n]_qA_{n-1}(q)$, where, for
instance, $q^{n-1}$ in the sum corresponds to having $1$ in the
rightmost position. Obviously, the set of these permutations is a
disjoint union of the set of all $n$-permutations that avoid $p$ and
the set of all $n$-permutations that quasi-avoid $p$. Thus,
$B^p_n(q)=[n]_qA^p_{n-1}(q)-A^p_n(q)$. Multiplying both sides of the
last equality by $x^n/[n]_q!$ and summing over all $n$ gives the
desired result.
\end{proof}

\begin{Lem}\label{lem02}\rm{(A $q$-analogue of~\cite[Thm. 25]{Kit2} that is valid for POPs)}
Let $P=p\mn\sigma$ be a POP, where $\sigma$ is an arbitrary POP
built on the alphabet that is incomparable to that involved in a
SPOP $p$. Then
$$A^P_q(x)=A^p_q(x)+A^{\sigma}_q(x)B^p_q(x).$$
\end{Lem}

\begin{proof} If a permutation $\pi$ avoids $p$ then it avoids $P$. Otherwise we
find the leftmost occurrence of $p$ in $\pi$. We assume that this
occurrence consists of the $|p|$ rightmost letters among the $i$
leftmost letters of $\pi$. So the subword of $\pi$ beginning at the
$(i+1)$st letter must avoid $\sigma$. From this, using independence
between the first $i$ letters of $\pi$ and the remain letters, we
conclude
$$A^P_n(q)=A^p_n(q)+\displaystyle\sum_{i=|\sigma|}^n
\left[\begin{array}{c} n \\ i \end{array}\right]_qB^p_i(q)
A^{\sigma}_{n-i}(q).$$ Observe that one can change the lower bound
in the sum above to 0, because $B^p_i(q)=0$ for
$i=0,1,\ldots,|p|-1$. Multiplying both sides by $x^n/[n]_q!$ and
summing over all $n$ we get the desired. \end{proof}

\begin{The}\label{mult}\rm{(A $q$-analogue of~\cite[Thm. 28]{Kit2} that is valid for POPs)} Let $p=p_1\mn\cdots\mn p_k$
be a multi-pattern ($p_i$s are SPOPs, and letters of $p_i$ and $p_j$
are incomparable for $i\neq j$). Then
$$A^p_q(x)=\displaystyle\sum_{i=1}^{k}A^{p_i}_q(x)\prod_{j=1}^{i-1}B^{p_j}_q(x)=
\sum_{i=1}^{k}A^{p_i}_q(x)\prod_{j=1}^{i-1}((x-1)A^{p_j}_q(x)+1).$$
\end{The}

\begin{proof} The first equality follows from lemma~\ref{lem02} by induction on
$k$, and the second equality is then given by
lemma~\ref{lem01}.\end{proof}

\begin{The}\rm{(A $q$-analogue of~\cite[Thm. 16]{Kit3})}\label{q-analog} If $N_p(\pi)$ denotes
the maximum number of non-overlapping occurrences of a SPOP $p$ in
$\pi$, then
$$\displaystyle\sum_{\pi}y^{N(\pi)}q^{inv(\pi)}\frac{x^{|\pi|}}{|\pi|!}=
\frac{A^p_q(x)}{1-yB^p_q(x)}=\frac{A^p_q(x)}{1-y((x-1)A^p_q(x)+1)}.$$\end{The}

\begin{proof} We fix $k$ and consider the multi-pattern $P_k=p\mn\cdots\mn p$ with
$k$ copies of $p$. A permutation avoiding $P_k$ has at most $k-1$
non-overlapping occurrences of $p$. From Theorem~\ref{mult},
$$A^{P_{k+1}}_q(x)-A^{P_k}_q(x)=A^p_q(x)(B^p_q(x))^k,$$ which is
a bivariate generating function for the number of permutations with
exactly $k$ non-overlapping occurrences of $p$ and with registered
inversions. The result is now follows from summing over all $k$ and
applying lemma~\ref{lem01}. \end{proof}

\section{Some open problems on POPs}\label{sec3}
We know very little on avoiding, and almost nothing on the
distribution of, POPs. There are a lot of posets and different
classes of posets, which provides enormous possibilities for further
research. In particular, a natural step would be to
extend/generalize results in the literature related to GPs to that
related to POPs in the manner Proposition~\ref{claes} and
Theorem~\ref{dist} are obtained. In this section, we state just few
problems on POPs that might be interesting to solve.

\subsection{Alternating patterns}

A permutation $\pi_1\pi_2\ldots\pi_n$ is {\em alternating} (resp.,
{\em reverse alternating}) if $\pi_1>\pi_2<\pi_3>\cdots$ (resp.,
$\pi_1<\pi_2>\pi_3<\cdots$). It is well known that the EGF for the
number of (reverse) alternating permutations is $\tan x+\sec x$.

We say that a permutation is a {\em $k$-non-alternating} (resp.,
{\em $k$-non-reverse-alternating}) if it does not contain $k$
consecutive letters that form an (resp., reverse) alternating
permutation. Using the complement, one can see that the numbers of
$k$-non-alternating and $k$-non-reverse-alternating $n$-permutations
are the same.

\begin{Prob} Enumerate $k$-non-alternating $n$-permutations. (For $k=4$ and $n\geq 4$ the numbers of ``good''
$n$-permutations are 19, 70, 331, 1863, 11637, 81110, ...; for $k=5$
and $n\geq 5$ we have the sequence 104, 528, 3296, 23168,
179712,...)
\end{Prob}

\begin{Prob} Enumerate $n$-permutations that are both
$k$-non-alternating and $k$-non-reverse-alternating. (For $k=4$ and
$n\geq 4$ we have the sequence 14, 52, 204, 1010, 5466, 34090,...;
for $k=5$ and $n\geq 5$ we have 24, 88, 458, 2716, 17808,
135182,...)\end{Prob}

\begin{figure}[h]
\begin{center}
\begin{picture}(10,3)
\put(-5,0){\put(12,2){\p} \put(14,0){\p} \put(16,2){\p}
\put(18,0){\p} \put(20,2){\p} \path(12,2)(14,0)(16,2)(18,0)(20,2)

\put(-1.2,0){$a_1$} \put(0.6,2){$a_2$} \put(3.6,0.8){$a_3$}
\put(6.3,2){$a_4$} \put(8.3,0){$a_5$}

\put(10.8,2){$a_1$} \put(13.6,0.8){$a_2$} \put(16.3,2){$a_3$}
\put(18.5,0){$a_4$} \put(20.3,2){$a_5$}

\put(0,0){\p} \put(2,2){\p} \put(4,0){\p} \put(6,2){\p}
\put(8,0){\p} \path(0,0)(2,2)(4,0)(6,2)(8,0)
 }
\end{picture}
\caption{Posets for the 5-reverse-alternating and 5-alternating
patterns.} \label{poset000}
\end{center}
\end{figure}
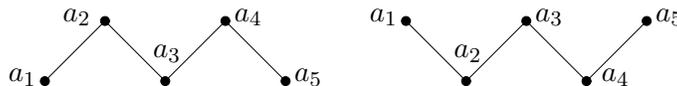

To generalize the problems above, we define a $k$-{\em alternating}
(resp., $k$-reverse-alternating) pattern to be one that forms a
(resp., reverse) alternating permutation of length $k$. Clearly, a
$k$-alternating (resp., $k$-reverse-alternating) segmented pattern
is a SPOP, where the corresponding poset is built on $k$ elements
$a_1,\ldots, a_k$ with the relations $a_1>a_2<a_3>\cdots$ (resp.,
$a_1<a_2>a_3<\cdots$) (see Figure~\ref{poset000} for the case
$k=5$).

Note that an occurrence of a descent in a permutation is an
occurrence of a 2-alternating pattern. Thus we have yet another
generalization of the notion of a descent beyond that discussed in
subsection~\ref{unimod}. Moreover, such patterns generalize the
patterns associated with peaks (valleys) in permutations, which
gives a motivation to study them.

The number of descents in a permutation $\pi$ is denoted by
$des(\pi)$. {\em Eulerian numbers} $A(n,k)$ count permutations in
the symmetric group $\SS_n$ with $k$ descents and they are the
coefficients of the {\em Eulerian polynomials} $A_n(t)$ defined by
$A_n(t)=\sum_{\pi\in\SS_n}t^{1+des(\pi)}$. The Eulerian polynomials
satisfy the identity
$$\sum_{k\geq 0}k^nt^k=\frac{A_n(t)}{(1-t)^{n+1}}.$$
For more properties of the Eulerian polynomials see~\cite{comtet}.

A natural generalization of the polynomials $A_n(t)$ is given by
considering $k$-alternating patterns instead of descents in the
definition of the polynomials. Let us call such new polynomials
$B^k_n(t)$. From definitions, $A_n(t)=B^2_n(t)$.

\begin{Prob} Study the properties of the polynomials $B^2_n(t)$ and
find the distribution for $k$-alternating patterns, that is, find an
explicit formula for $B^2_n(t)$ or, if possible, coefficients of
$B^2_n(t)$.
\end{Prob}

\begin{Prob} Find joint distribution for $k$-alternating and $k$-reverse-alternating patterns. \end{Prob}

Note that there are many other (segmented) patterns that can be
built on posets similar to those in Figure~\ref{poset000}. For
example, one could consider the pattern $a_1a_3a_5a_2a_4$ built on
the poset to the left in Figure~\ref{poset000}. To study other than
alternating patterns built on such posets might be also an
interesting direction to explore.

\subsection{Co-unimodal patterns}

Recall from subsection~\ref{unimod} that a SPOP
$\sigma=\sigma_1\sigma_2\ldots\sigma_k$ is co-unimodal if
$k=\sigma_1>\sigma_2>\cdots
>\sigma_j<\cdots <\sigma_k$ for some $2\leq j\leq k$. We extend
the concept of co-unimodal pattern to that of {\em free co-unimodal
pattern} by removing the restriction ``$k=$'' in the definition.
Note that co-unimodal patterns impose weaker restrictions on
permutations than free co-unimodal patterns do.

\begin{Prob} How many of $n$-permutations avoid a co-unimodal pattern of length~$k$.
(For $k=4$ and $j=2$ (resp., $j=3$) see the record for the pattern
$utxv$ (resp., $spor$) in table~\ref{unsolved}.)\end{Prob}

\begin{Prob} How many of $n$-permutations avoid a free co-unimodal pattern of length $k$.
(For $k=4$, because of the complement, $j=2$ and $j=3$ give the same
number of $n$-permutations avoiding them; see the record for the
pattern $ijkm$ in table~\ref{unsolved}.)\end{Prob}

\begin{Prob} Find the distribution of a co-unimodal pattern of length $k$. \end{Prob}

\begin{Prob} Find the distribution of a free co-unimodal
pattern of length~$k$.\end{Prob}

\begin{Prob} Find the number of $n$-permutations avoiding simultaneously two or more of (free)
co-unimodal patterns. We provide some numerical data in case $k=4$.
Suppose $F_2$ and $F_3$ are the free co-unimodal patterns
corresponding to $j=2$ and $j=3$ respectively; also, $U_2$ and $U_3$
are the co-unimodal patterns corresponding to $j=2$ and $j=3$
respectively. The initial values for the number of $n$-permutations,
$n\geq 4$, avoiding a pair of the patterns are as follows:
($F_2$,$F_3$) -- 18, 66, 252, 1176, 5768, 34216; ($F_2$,$U_3$) --
19, 75, 330, 1753, 10319, 70011; ($F_3$,$U_2$) -- 20, 81, 372, 1981,
11866, 80043; ($U_2$,$U_3$) -- 21, 91, 462, 2718, 18181, 136491.
\end{Prob}

\begin{Prob} Find the joint distribution of two or more (free) co-unimodal patterns. \end{Prob}

\subsection{Remaining cases of SPOPs of length
four}\label{usolved-4}

\begin{figure}[h]
\begin{center}
\begin{picture}(6,4)
\put(0,0){\put(2,0){\p} \put(2,2){\p} \put(4,0){\p} \put(4,2){\p}
\put(4,4){\p} \path(2,0)(2,2) \path(4,0)(4,2)(4,4)

\put(1,0){$1'$} \put(1,2){$2'$} \put(4.3,0){1} \put(4.3,2){2}
\put(4.3,4){3}
 }
\end{picture}
\caption{Poset from which some patterns in table~\ref{unsolved} are
built.} \label{poset001}
\end{center}
\end{figure}
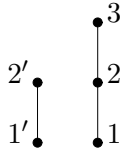

In table~\ref{unsolved}, we record few initial values for the number
of $n$-permutations in some of unsolved cases of avoidance of SPOPs
of length four, $n\geq 1$. In the table we record patterns having at
least one pair of incomparable letters (see Figure~\ref{poset001}
for the corresponding poset), although there are unsolved cases when
all elements are comparable (we have a chain in the Hasse diagram).
We refer to~\cite{ElizNoy} for information on unsolved segmented GPs
of length four. Table~\ref{unsolved} is also an extended version of
the corresponding table in~\cite{Kit3}.

\begin{table}[ht]
\begin{center}
   \begin{tabular}{|l|l|}
    \hline
$11'22'$ & $1,\ 2,\ 6,\ 18,\ 70,\ 300,\ 1435,\ 7910,\
47376,\ldots$\\
   \hline
$121'2'$ & $1,\ 2,\ 6,\ 18,\ 61,\ 281,\ 1541,\ 8920,\
57924,\ldots$\\
   \hline
$11'2'2$ & $1,\ 2,\ 6,\ 18,\ 71,\ 322,\ 1665,\ 9789,\
64327,\ldots$ \\
   \hline
$12'1'2$ & $1,\ 2,\ 6,\ 18,\ 61,\ 272,\ 1410,\ 8048,\
51550,\ldots$\\
   \hline
$121'3$ & $1,\ 2,\ 6,\ 20,\ 83,\ 411,\ 2290,\ 14588,\
104448,\ldots$\\
   \hline
$131'2$ & $1,\ 2,\ 6,\ 20,\ 81,\ 390,\ 2161,\ 13678,\
96983,\ldots$\\
   \hline
$231'1$ & $1,\ 2,\ 6,\ 20,\ 83,\ 402,\ 2245,\ 14192,\
100650,\ldots$\\
   \hline
$abcd$ & $1,\ 2,\ 6,\ 19,\ 70,\ 331,\ 1863,\ 11637,\
81110,\ldots$\\
   \hline
$utxv$ & $1,\ 2,\ 6,\ 23,\ 110,\ 630,\ 4210,\ 32150,\
276210,\ldots$\\
   \hline
$spor$ & $1,\ 2,\ 6,\ 22,\ 100,\ 540,\ 3388,\ 24248,\
195048,\ldots$\\
   \hline
$ijkm$ & $1,\ 2,\ 6,\ 21,\ 90,\ 450,\ 2619,\ 17334,\
129114,\ldots$\\
   \hline
$egfh$ & $1,\ 2,\ 6,\ 20,\ 84,\ 412,\ 2300,\ 14676,\
104536,\ldots$\\
   \hline
$efgh$ & $1,\ 2,\ 6,\ 20,\ 80,\ 404,\ 2368,\ 15488,\
114480,\ldots$\\
   \hline
$fegh$ & $1,\ 2,\ 6,\ 20,\ 80,\ 360,\ 1888,\ 11168,\
75168,\ldots$\\
   \hline
   \end{tabular}
\smallskip \caption{The initial values  for the number of $n$-permutations
avoiding 4-SPOPs in a few of unsolved cases, $n\geq 1$. See
Figures~\ref{poset001} and~\ref{poset010} for the corresponding
poset.} \label{unsolved}
\end{center}
\end{table}

\begin{figure}[h]
\begin{center}
\begin{picture}(2,4)
\put(-10,0){\put(5,0){\p} \put(5,2){\p} \put(7,0){\p} \put(7,2){\p}
\path(5,0)(5,2)(7,0)(7,2)(5,0)

\put(-0.8,0){$a$} \put(-0.7,2){$b$} \put(2.3,0){$c$}
\put(2.3,2){$d$}

\put(4.3,0){$e$} \put(4.3,2){$f$} \put(7.3,0){$g$} \put(7.3,2){$h$}

\put(9.3,2){$i$} \put(11.3,1.2){$j$} \put(12.7,0.3){$k$}
\put(15.3,2){$m$}

\put(9,2){\p} \put(11,1){\p} \put(13,0){\p} \put(15,2){\p}
\path(9,2)(11,1)(13,0)(15,2)

\put(0,0){\p} \put(0,2){\p} \put(2,0){\p} \put(2,2){\p}
\path(0,0)(0,2)(2,0)(2,2)

\put(16.3,0.5){$p$} \put(17.1,-0.2){$o$} \put(17.3,2){$s$}
\put(19.3,1){$r$}

\put(18,0){\p} \put(17,1){\p} \put(18,2){\p} \put(19,1){\p}
\path(18,0)(17,1)(18,2)(19,1)(18,0)

\put(21,0){\p} \put(21,2){\p} \put(23,1.5){\p} \put(23,0.5){\p}
\path(21,0)(21,2)(23,1.4)(23,0.6)(21,0)

\put(20.3,0){$t$} \put(20.2,2){$u$} \put(23.3,1.7){$v$}
\put(23.3,0.4){$x$}

 }
\end{picture}
\caption{Five posets to build 4-SPOPs that were not considered.}
\label{poset010}
\end{center}
\end{figure}
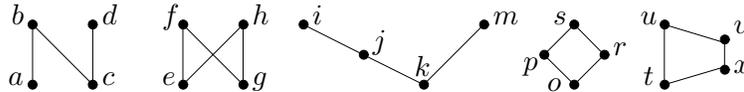

Other 4-SPOPs that were not considered can be built on the posets in
Figure~\ref{poset010}. For example, for the second poset there are
three SPOPs to consider that are non-equivalent up to trivial
bijections: $egfh$, $efgh$, and $fegh$ (see table~\ref{unsolved} for
corresponding sequences).

Notice that the leftmost poset in Figure~\ref{poset010} can be used
to build the 4-reverse-alternating pattern $abcd$, as well as the
4-alternating pattern $dcba$, whereas the third (resp., forth,
fifth) poset in Figure~\ref{poset010} can be used to build free
co-unimodal (resp., co-unimodal) pattern(s) of length 4, namely
$ijkm$ (resp., $spor$, $utxv$).

\subsection{Further research directions}

The problems stated above can be extended to many POPs by inserting
dash(es) in the SPOPs discussed. Also, a natural generalization of
any avoidance problem is finding the distribution of a (S)POP under
consideration. Moreover, joint distribution of (S)POPs and, as a
special case, multi-avoidance of these patterns, is a possible
direction for further research after choosing (S)POPs to consider.
All these problems are interesting from enumerative point of view
but also might bring interesting connections to other combinatorial
objects, in which case, as always, explicit bijections would be
desirable.

\section{Acknowledgements}

The author is grateful to Ira Gessel for the discussion and
references related to the distribution of peaks in permutations; to
Jeff Remmel for the discussion on $q$-analogues for non-overlapping
occurrences of patterns, as well as for his support during the
author's stay at UCSD; to the two anonymous referees for their
helpful comments.

\end{document}